\documentclass[]{elsarticle}
\usepackage{amssymb,amsmath,amsthm,amsfonts}
\usepackage{verbatim}

\usepackage{color}  
\usepackage[colorlinks=true, linkcolor=blue, citecolor=blue]{hyperref}  

\newtheorem{theorem}{Theorem}[section]  
\newtheorem{lemma}{Lemma}[section]
\newtheorem{corollary}{Corollary}[section]
\newtheorem{proposition}{Proposition}[section]
\newtheorem{problem}{Problem}[section]
\newtheorem{example}{Example}[section]
\newtheorem{definition}{Definition}[section]
\newtheorem{remark}{Remark}[section]

\newcommand{\thm}[2]{\begin{theorem}\label{#1}{\sl #2}\end{theorem}}
\newcommand{\lem}[2]{\begin{lemma}\label{#1}{\sl #2}\end{lemma}}
\newcommand{\cor}[2]{\begin{corollary}\label{#1}{\sl #2}\end{corollary}}

\newcommand{\rem}[2]{\begin{remark}\label{#1}{\rm #2}\end{remark}}

\newcommand{\lt}{\left}
\newcommand{\rt}{\right}
\newcommand{\toi}{\to\infty}
\newcommand{\oF}{\overline{F}}

\newcommand{\be}{\begin{eqnarray*}}
\newcommand{\ee}{\end{eqnarray*}}
\newcommand{\pf}{\noindent\textbf{Proof. }}

\begin{document}

\begin{center}
{\bf \Large
{Uniform asymptotics for the tail probability of weighted sums with heavy tails}}
\end{center}

\centerline{Chenhua Zhang
        \footnote{Department of Mathematics, The University of Southern Mississippi, Hattiesburg, MS 39406-5045, USA, chenhua.zhang@usm.edu}
     }

{\narrower\small
\baselineskip=11pt\parindent=0pt

\bigskip
 {\bf Abstract.}
This paper studies the tail probability of weighted sums of the form $\sum_{i=1}^n c_i X_i$, where random variables $X_i$'s are either independent or pairwise quasi-asymptotical independent with heavy tails.
Using $h$-insensitive function, the uniform  asymptotic  equivalence of the tail probabilities of $\sum_{i=1}^n c_iX_i$, $\max_{1\le k\le n}\sum_{i=1}^k c_iX_i$ and $\sum_{i=1}^n c_iX_i^+$ is established, where $X_i$'s  are independent and follow the long-tailed distribution, and $c_i$'s take value in a broad interval. Some further uniform asymptotic results for the weighted sums of $X_i$'s with dominated varying tails are obtained. An application to the ruin probability in a discrete-time insurance risk model is presented.

\bigskip


\bigskip

}

\bigskip
\section{Introduction}
In this paper, all asymptotic and limit relations are taken as $x\toi$ unless otherwise stated.
For independently and identically distributed (iid) subexponential random variables $X_i,i\ge1$, it is well-known that, for any $n\ge2$,
\begin{eqnarray}\label{S_n sim M_n}
 P\bigg(\sum_{i=1}^n   X_i>x\bigg)\sim P\bigg(\max_{1\le k\le n}\sum_{i=1}^k   X_i>x\bigg)  \sim P\bigg(\sum_{i=1}^n   X_i^+>x\bigg)\sim  \sum_{i=1}^n P(
 X_i>x),
\end{eqnarray}
where $x^+=\max\{x,0\}$. There are quite a few ways to generalize these asymptotic relations. One way is to consider some broader classes of
heavy-tailed distributions, see, e.g., Ng et al. \hyperref[]{\cite{NTY02}}. Another way is to study the randomly stopped sums, see, e.g., Denisov et al.
\hyperref[]{\cite{DFK10}}. Allowing some dependence of $X_i$'s, similar results can be obtained for different classes of heavy-tailed distributions, see
Wang and Tang \hyperref[]{\cite{WT04}},  Geluk and Ng \hyperref[]{\cite{GN06}}, Tang \hyperref[]{\cite{T08}} , Geluk and Tang \hyperref[]{\cite{GT09}},
and references therein.

A more general way is to work on the weighted sums of form $\sum_{i=1}^n c_iX_i$, where weights $c_i$'s are real numbers. If $X_i$'s are iid
subexponential random variables, Tang and Tsitsiashvili \hyperref[]{\cite{TT03}} proved that for any $0<a\le b<\infty$, the asymptotic relation
\begin{eqnarray}\label{expansion for weighted sum of iid in S}
 P\bigg(\sum_{i=1}^n c_i X_i>x\bigg)\sim  \sum_{i=1}^n P(c_iX_i>x),
\end{eqnarray}
holds uniformly for  $a\le c_i\le b, 1\le i\le n,$ in the sense that
\be
  \lim_{x\toi}\sup_{a\le c_i\le b, 1\le i\le n}\left|\frac{P(\sum_{i=1}^n c_i X_i>x)}{\sum_{i=1}^n P(c_iX_i>x)}-1\right|=0.
\ee
Recently, Liu et al. \hyperref[]{\cite{LGW12}} and Li \hyperref[]{\cite{L13}} established the same asymptotic relation for some dependent $X_i$'s.

Chen et al.  \hyperref[]{\cite{CNY11}} showed that for any fixed $0<a\le b<\infty$ it holds that uniformly for $a\le c_i\le b$, $1\le i\le n$,
\begin{eqnarray}\label{S_n sim S_{(n)} sim S_n^+}
   P\bigg(\sum_{i=1}^n c_iX_i>x\bigg)\sim P\bigg(\max_{1\le k\le n}\sum_{i=1}^k c_iX_i>x\bigg) \sim P\bigg(\sum_{i=1}^n c_iX_i^+>x\bigg),
\end{eqnarray}
where $X_i$'s are independent, not necessarily identically distributed, random variables with long-tailed distributions.
This result is extended by substituting $b$ with any positive function $b(x)$ such that $h(x)\nearrow\infty$ and $b(x)=o(x)$  in this paper.

Replacing the constant weights $c_i$'s with random weights $\theta_i$'s,  the asymptotic relation \hyperref[]{\eqref{expansion for weighted sum of iid
in S}} and \hyperref[]{\eqref{S_n sim S_{(n)} sim S_n^+}} still hold if the weights $\theta_i$'s, independent of $X_i$'s, are uniformly bounded away from
zero and infinity. Then it is very natural to consider the randomly weighted sum of form $\sum_{i=1}^n \theta_iX_i$.
Wang and Tang \hyperref[]{\cite{WT06}} obtained
$   P\big(\sum_{i=1}^n \theta_iX_i>x\big)\sim P\big(\max_{1\le k\le n}\sum_{i=1}^k \theta_iX_i>x\big) \sim P\big(\sum_{i=1}^n \theta_iX_i^+>x\big)
$
for the case that the random weights are not necessarily bounded and $X_i$'s are independently random variables with common distribution belonging to a
smaller class than the class of subexponential distributions. Furthermore, Zhang et al. \hyperref[]{\cite{ZSW09}}, Chen and Yuen \hyperref[]{\cite{CY09}}
established the same results for dependent $X_i$'s, where the  dependence structures of $X_i$'s are essentially same for proof of their results.

The rest of this paper is organized as follows. Section 2 reviews some important classes of heavy-tailed distributions. Section 3 states the main results
along with some corollaries. Section 4 gives an application of the main results to the ruin probability in a discrete-time insurance risk model. The proof
of the main results and some lemmas are presented in Section 5.

\section{Classes of Heavy-Tailed Distributions}

A random variable $X$ or its distribution $F$ is said to be heavy-tailed to the right or have a heavy (right) tail if the corresponding moment generate function does not exist on the
positive real line, i.e., $E e^{tX}=\int_{-\infty}^\infty e^{tx}d F(x)=\infty$ for any $t>0$.
The most important class of heavy-tailed distributions is the class of subexponential distributions, denoted by $\mathcal{S}$. Write the tail
distribution by $\oF(x)=1-F(x)$ for any distribution  $F$. Let $F^{*n}$ denote the $n$-fold convolution of $F$.
A distribution $F$ concentrated on $[0,\infty)$ is subexponential if
\be
\overline{ F^{*n}}(x)\sim n\overline{F}(x)
\ee
 for some or, equivalently, for all $n\ge2$.
 More generally, a distribution $F$ on $(-\infty,\infty)$ belongs to the subexponential class if $F^+(x)=F(x)I_{\{x\ge0\}}$ does.

 Closely related to the subexponential class $\mathcal{S}$,  the class $\mathcal{D}$ of dominated varying distributions consists of distributions satisfying
 \be
 \limsup_{x\toi}\frac{\oF(yx)}{\oF(x)}<\infty
\ee
  for some or, equivalently, for all $0<y<1$.
A slightly smaller class of $\mathcal{D}$ is the class of distributions with consistently varying tail, denoted by $\mathcal{C}$. Say that a distribution
$F$ belongs to the class $\mathcal{C}$ if
\be
\lim_{y \searrow1}\liminf_{x\toi}\frac{\overline{F}(yx)}{\overline{F}(x)}=1 \mathrm{\ or,\ equivalently,\ } \lim_{y
\nearrow1}\limsup_{x\toi}\frac{\overline{F}(yx)}{\overline{F}(x)}=1.
\ee

A distribution $F$ belongs to the class $\mathcal{L}$ of long-tailed distributions if
 \be
 \lim_{x\toi}\frac{\oF(x+y)}{\oF(x)}=1
\ee
for some or, equivalently, for all $y$. A tail distribution $\oF$ is called $h$-insensitive if $\oF(x+ y)\sim\oF(x)$ holds uniformly for all $|y|\le
h(x)$, where $h(x)$ is a positive nondecreasing function and $\lim_{x\toi}h(x)=\infty$.
The concept of $h$-insensitive function is extensively used in the monograph of Foss et al. \hyperref[]{\cite{FKZ11}}.
For any distribution $F\in\mathcal{L}$, it can be shown that $\oF$ is $h$-insensitive for some positive nondecreasing function $h(x):=h_F(x)$ such that
$h(x)\nearrow\infty$ and  $h(x)=o(x)$, see, e.g., Lemma \hyperref[]{\ref{h-insensitive for F in L}} in Section \ref{Proof Section},  Section 2 in Foss and Zachary \hyperref[]{\cite{FZ03}}, Lemma 4.1 of Li et al. \hyperref[]{\cite{LTW10}}. Consequently,
 $\oF$ is $ch$-insensitive for any fixed positive real number $c$.

It is known that  the proper inclusion relations $$\mathcal{C} \subset \mathcal{D} \cap \mathcal{L} \subset  \mathcal{S}\subset  \mathcal{L}$$ hold,
 see, e.g., Embrechts et al. \hyperref[]{\cite{EKM97}}, Foss et al. \hyperref[]{\cite{FKZ11}}.

\section{Main Results}

Throughout the rest of this paper $X_i, i\ge 1$, are random variables with distribution $F_i, i\ge1$, respectively. Adopt the notation $M_{c}F$ and $*_{1\le i\le n}M_{c_i}F_i$ in Barbe  and McCormick \hyperref[]{\cite{BM09}}. For $X\sim F$ and $c>0$, let $M_{c}F(x)= F(x/c)$ be the distribution of $cX$. The distribution of $\sum_{i=1}^nc_iX_i$ is $*_{1\le i\le n}M_{c_i}F_i$, where $X_i,1\le i\le n,$ are independent random variables and $*_{1\le i\le n}M_{c_i}F_i$ is the convolution of $ M_{c_i}F_i,1\le i\le n$.

The first main result generalizes Lemma 4.1 of Chen et al. \hyperref[]{\cite{CNY11}} with different approach in two ways. First, it increases the upper  bound of the weights and decreases the lower bound of the weights. Second, the fixed shift term $A$ in Lemma 4.1 of Chen et al. \hyperref[]{\cite{CNY11}} is enlarged to some unbounded function, which is irrespective of the upper bound of the weights.
\thm{}{\label{uniform h-insensitive for convolution in L}
If $X_i\sim F_i \in \mathcal{L},1\le i\le n$, are independent random variables,
there exists a positive nondecreasing function $h(x):=h(x;F_1,\cdots,F_n)$ satisfying $h(x)\nearrow\infty$ such that $*_{1\le i\le n}M_{c_i}F_i$ is uniformly $h(x)$-long-tailed for $a(x)\le c_i\le  b(x), 1\le i\le n$, in the sense that
\be
P\bigg(\sum_{i=1}^nc_iX_i>x\pm h(x)\bigg)\sim P\bigg(\sum_{i=1}^nc_iX_i>x \bigg)
\ee
holds uniformly for $a(x)\le c_i\le  b(x), 1\le i\le n$, i.e.,
\begin{eqnarray}\label{uniform h-insensitive}
 &&\lim_{x\toi} \sup_{a(x)\le c_i\le  b(x), 1\le i\le n}
        \left|\frac{\overline{*_{1\le i\le n}M_{c_i}F_i}(x\pm h(x))}{\overline{*_{1\le i\le n}M_{c_i}F_i}(x)}-1\right| =0,
\end{eqnarray}
where  the positive function $b(x)$  satisfies $b(x)\nearrow\infty$ and $b(x) = o(x)$, $h(x)$ is irrespective of $b(x)$, $a(x)=h^{-\delta}(x)\searrow0$ for some $\delta>0$.
}

\rem{}{ Considering the case of Weibull distribution $F_1(x)=1-e^{-cx^{\tau}}\in\mathcal{S}\subset\mathcal{L}$ with $0<\tau<1$, it  indicates that the restriction on $a(x)$  can not be weakened in general.
}
It is known that the class  $\mathcal{L}$ is closed under convolution (see, e.g., Theorem 3 of Embrechts and Goldie \hyperref[]{\cite{EG80}}, Corollary
2.42 of Foss et al. \hyperref[]{\cite{FKZ11}}), which can be also derived directly from Theorem \ref{uniform h-insensitive for convolution in L}.

\cor{}{
If $X_i\sim F_i\in\mathcal{L}, 1\le i\le n$, are independent random variables, then the distribution of $\sum_{i=1}^nc_iX_i>x \big)$  is long-tailed for any fixed $c_i>0, 1\le i\le n$. Consequently, the class $\mathcal{L}$ of long-tailed distributions is closed under convolution.}

\thm{}{\label{2nd theorem}
If $X_i\sim F_i \in \mathcal{L},1\le i\le n$, are independent random variables,
there exist positive functions $a(x)$ and $b(x)$ satisfying $a(x)\searrow0$ and $b(x)\nearrow\infty$ such that  the asymptotic relations \hyperref[]{\eqref{S_n sim S_{(n)} sim S_n^+}} hold uniformly for $a(x)\le c_i\le b(x)$, $1\le i\le n$.
}
  The following result can be also founded in  Lemma 3.4 of Foss et al. \hyperref[]{\cite{FKZ11}}.
\cor{}{\label{S in R}
A distribution $F\in\mathcal{S}$ iff $F\in\mathcal{L}$ and $\overline{F*F}(x)\sim 2\oF(x)$.
}

  Random variables $X_i,i\ge1,$ are  pairwise strong quasi-asymptotically independent (pSQAI)
if, for any $i\ne j$,
\be
\lim_{\min\{x_i,x_j\}\to\infty}P\lt( |X_i|> x_i|X_j> x_j\rt)=0,
\ee
which was used in Geluk and Tang \hyperref[]{\cite{GT09}}, Liu et al. \hyperref[]{\cite{LGW12}} and Li \hyperref[]{\cite{L13}},   and related to what is called asymptotic independence; see e.g.  Resnick \hyperref[]{\cite{R02}}.

\thm{}{ \label{3rd theorem}
If $X_i\sim F_i\in\mathcal{C},1\le i\le n$, are pSQAI random variables  and $b(x)$ is an arbitrary fixed positive function
satisfying $b(x)\nearrow\infty$ and $b(x)=o(x)$,
then it holds  that, uniformly for any  $0<c_i\le b(x), 1\le i\le n,$
\begin{eqnarray}\label{Tail equivalance}
   P\bigg(\sum_{i=1}^n c_iX_i>x\bigg)\sim P\bigg(\max_{1\le k\le n}\sum_{i=1}^k c_iX_i>x\bigg) \sim P\bigg(\sum_{i=1}^n c_iX_i^+>x\bigg)\sim \sum_{i=1}^n
   P(c_iX_i>x).
\end{eqnarray}
}
\cor{}{\label{Cor to 3rd Theorem}
Under assumption of Theorem \ref{3rd theorem}, the above result still holds for $0\le c_i\le b(x), 1\le i\le n,$ and $\min_{1\le i\le n}c_i>0$.
}

The next theorem extends  Lemma 2.1 of Liu et al \hyperref[]{\cite{LGW12}} and Theorem 2.1 of Li \hyperref[]{\cite{L13}} with a different proof, which is based on Theorem \ref{uniform h-insensitive for convolution in L}.

\thm{}{ \label{4th theorem}
If  $X_i\sim F_i\in\mathcal{D}\cap\mathcal{L},1\le i\le n,$ are pSQAI random variables, there exist  a positive function $a(x)\searrow0$ and a positive function $b(x)\nearrow\infty$ 
such that \hyperref[]{\eqref{Tail equivalance}} holds uniformly for $a(x)\le c_i\le b(x), 1\le i\le n$.
}

\rem{}{ Both $a(x)$ and $b(x)$ depend on $h(x)$ in Theorem \ref{2nd theorem} and \ref{4th theorem}, where $h(x)=o(x)$ is given in Theorem \ref{uniform h-insensitive for convolution in L}. More specifically, $a(x)= h^{-\delta}(x)$ for some $\delta>0$ and  $b(x)=o(h(x))$, for example, $b(x)=h^{1/2}(x)$.
}

\rem{}{ If the constant weights $c_i,1\le i\le n$ are replaced by random weights $\theta_i,1\le i\le n$, which are independent of $X_i,1\le i\le n$,
conditioning on the random weights can easily establish the corresponding results for random weights sums.
}

The proof of Theorem \ref{4th theorem} gives an extension of Lemma 4.3 of Geluk and Tang \hyperref[]{\cite{GT09}}.
\cor{}{\label{liminf for L}
If $X_i\sim F_i \in \mathcal{L},1\le i\le n$, are pQSAI random variables, 
 it holds that, for some the positive functions $b(x)\nearrow\infty$ and  $a(x)\searrow0$,
\begin{eqnarray}\label{uniform h-insensitive}
 \lim_{x\toi} \inf_{a(x)\le c_i\le  b(x), 1\le i\le n}
        \frac{P\big(\sum_{i=1}^nc_iX_i>x\big)}{\sum_{i=1}^nP\big(c_iX_i> x\big)}\ge 1.
\end{eqnarray}
}
\section{Application to Risk Theory}
Consider the following discrete-time insurance risk model
\be
U_0=x,\ U_n=U_{n-1}(1+r_n)-X_n,n\ge 1,
\ee
where $U_n$ stands an insurer's surplus at the end of period $n$ with a deterministic initial surplus $x$, $r_n$ represents the constant interest force
of an insurer's risk-free investment, and
the net loss $X_n$ over period $n$  equals the total amount of claims plus other costs minus the total amount of premiums during period $n$.
It is an interesting and important problem arising from the above discrete-time insurance risk model  to study the ruin probabilities of the insurer.
See Tang \hyperref[]{\cite{T04}} for detailed discussion.

The ruin probability by time $n$ is defined as
\be
\psi(x;n)=P\Big(\min_{i=1}^{n} U_i<0\,|\,U_0=x\Big).
\ee
It is easy to see that the  surplus process is of form
\be
U_0=x,\ U_n=\prod_{i=1}^n (1+r_i) x -\sum_{i=1}^n \Big(\prod_{j=i+1}^n (1+r_j) \Big) X_i, n\ge1.
\ee
Define the discounted surplus process as follows
\be
\widetilde{U}_n=\Big(\prod_{i=1}^n (1+r_i)\Big)^{-1}U_n=  x -\sum_{i=1}^n c_i X_i,
\ee
where $c_i=\prod_{j=1}^i (1+r_j)^{-1}$ represents the discount factor from time $i$ to time $0$, $1\le i\le n$.
Then the corresponding ruin probability can be written as
\be
\psi(x;n)=P\Big(\min_{i=1}^{n} \widetilde{U}_i<0\,|\,\widetilde{U}_0=x\Big)
= P\Big(\max_{1\le i\le k}\sum_{i=1}^k c_i X_i>x\Big).
\ee
 Applying  Theorem \hyperref[]{\ref{2nd theorem}} and Theorem \hyperref[]{\ref{4th theorem}} in Section 3, the
following asymptotic results can be obtained.
\cor{}{
Assume that net losses $X_i, i\ge 1$ are independent random variables, which are not necessarily identically distributed, with distribution $F_i,i\ge1$,
respectively.
If  $F_i\in\mathcal{L},1\le i\le n$, then
\begin{eqnarray*}
   \psi(x;n) \sim P\bigg(\sum_{i=1}^n c_iX_i>x\bigg)\sim  P\bigg(\sum_{i=1}^n c_iX_i^+>x\bigg).
\end{eqnarray*}
If  $F_i\in\mathcal{D}\cap\mathcal{L},1\le i\le n$, then
\begin{eqnarray*}
   \psi(x;n)\sim P\bigg(\sum_{i=1}^n c_iX_i>x\bigg)\sim  P\bigg(\sum_{i=1}^n c_iX_i^+>x\bigg)\sim \sum_{i=1}^n P(c_iX_i>x).
\end{eqnarray*}
}

\section{Proof of Results}\label{Proof Section}

A function $h(x)$ is called slowly varying at infinity if $h(xy)\sim h(x)$ for any $y>0$, It is well-known that $h(x)=o(x^{\delta})$ for any $\delta>0$ if $h(x)$ is a slowly varying function, see, e.g., Bingham et al. \hyperref[]{\cite{BGT89}}.
The following result is crucial for the proof of all theorems in this paper. It shows that any tail distribution of a long-tailed distribution is uniformly $h$-insensitive for a slowly varying function $h$.
\lem{}{\label{h-insensitive for F in L}
If $X\sim F\in\mathcal{L}$, then $\oF$ is $h$-insensitive for a positive nondecreasing and slowly varying function $h(x):=h(x;F): (0,\infty)\to(0,\infty)$ satisfying $h(x)\nearrow\infty$, $h(x)\le c h(\frac{x}{c})$ for all $c\ge 1$, and
\begin{eqnarray}\label{Theorem3.1.n=1}
 \lim_{x\toi} \sup_{a(x)\le c\le  b(x)} \left|\frac{P(cX>x\pm h(x))}{P(cX>x)}-1\right|=0,
\end{eqnarray}
where   $b(x)$ is an arbitrary positive function such that $b(x)\nearrow\infty$ and $b(x)=o(x)$, and $a(x)=h^{-\delta}(x)$ for some $\delta>0$.
}
\pf
For any fixed $\delta>0$, let $\{x_n,n\ge 1\}$ be a sequence of increasing positive real numbers such that $x_{n+1} \ge 2x_{n}>0$, $n\ge 1$, and for any $x\ge x_n$,
\begin{eqnarray}\label{step 1 to define h function}
\sup_{|y|\le n}\left|\frac{\overline{F}(x+y)}{\overline{F}(x)}-1\right|
\le  \max\left\{\left|\frac{\overline{F}(x+n^{1+\delta})}{\overline{F}(x)}-1\right|,  \left|\frac{\overline{F}(x-n^{1+\delta})}{\overline{F}(x)}-1\right|\right\} \le
\frac{1}{n}.
\end{eqnarray}
Borrowing the idea of the proof of Corollary 2.5 in \hyperref[]{\cite{CS94}}, let
$$
h(x)=\left\{ \begin{array}{ll}
 \frac{2}{x_1}x& x_0=0< x<x_1\\
 n+\frac{x-x_{n-1}}{x_n-x_{n-1}}&x_{n-1}\le x<x_{n},n\ge2.
 \end{array}\right.
$$
Clearly, $h(x)$ is a positive nondecreasing, piecewise linear, continuous function and $h(x)\nearrow\infty$.
Since $h(x)$ is a nondecreasing function, $h(xy)\sim h(x)$ for any $y>0$ is equivalent to $h(2x)\sim h(x)$, which follows from the facts that $h(x)\nearrow\infty$ and $h(x)\le h(2x)<h(x_{n+1})=n+2\le h(x)+2$ for any $x_{n-1}\le x<x_{n}$.

\noindent
For any $x\ge x_n$, i.e., $x\in [x_{n+k},x_{n+k+1})$ for some $k:=k(x)\ge0$, and $|y|\le h^{1+\delta}(x)=(n+k+1)^{1+\delta}$,
it follows from \hyperref[]{\eqref{step 1 to define h function}} that
\begin{eqnarray*}
\sup_{|y|\le h^{1+\delta}(x)}\left|\frac{\overline{F}(x+y)}{\overline{F}(x)}-1\right|
\le \frac{1}{ n+k+1 }
\le \frac{1}{n }\to 0,\quad \mathrm{as}\ n\toi,
\end{eqnarray*}
i.e., $\oF$ is $h^{1+\delta}$-insensitive, which of course implies that $\oF$ is $h$-insensitive.
 Since $x_{n+1}-x_n\ge x_n\ge x_{n}-x_{n-1}, n\ge 1$,  $h'(x)$ is a nonincreasing function on $ \cup_{n=1}^\infty (x_{n-1},x_{n})$, which implies that
 $h(x)$ is a concave function on $[0,\infty)$.  The concavity of $h(x)$ and the fact $h(0)=0$ lead to
 $h(\frac{x}{c})=h\big(\frac{1}{c}x+(1-\frac{1}{c})0\big)\ge \frac{1}{c}h(x)+(1-\frac{1}{c})h(0)=\frac{1}{c}h(x)$, i.e.,  $h(x)\le ch(\frac{x}{c})$, for
 any $x>0,c>1$.\\
Hence, $\frac{h(x)}{c}\le h\big(\frac{x}{c}\big)  \le h^{1+\delta}(\frac{x}{c}\big)$ for $1\le c\le b(x)$.
 Note that $\frac{h(x)}{c}\le\frac{h(x)}{a(x)} =h^{1+\delta}(x)\le h^{1+\delta}(\frac{x}{c}\big)$ for $a(x)\le c\le 1$.
 The monotonicity of $\oF$ yields $\overline{F}\big(\frac{x}{c}+  h^{1+\delta}(\frac{x}{c})\big)\le  P\big(cX>x\pm h(x)\big)=\overline{F}\big(\frac{x}{c}\pm \frac{h(x)}{c}\big) \le \overline{F}\big(\frac{x}{c}-h^{1+\delta}(\frac{x}{c})\big)$ for $a(x)\le c\le b(x)$.
The uniform asymptotic relation \hyperref[]{\eqref{Theorem3.1.n=1}}  follows from the inequalities
\be
\frac{\overline{F}\big(\frac{x}{c}+h^{1+\delta}(\frac{x}{c}\big)\big)}{\overline{F}\big(\frac{x}{c}\big)}-1
&\le& \frac{P(cX>x\pm h(x))}{P(cX>x)}-1=\frac{\overline{F}\big(\frac{x}{c}\pm \frac{h(x)}{c})\big) }{\overline{F}\big(\frac{x}{c}\big)}-1\\
&\le& \frac{\overline{F}\big(\frac{x}{c}-h^{1+\delta}(\frac{x}{c}\big)\big)}{\overline{F}\big(\frac{x}{c}\big)}-1,\quad a(x)\le c\le b(x),
\ee
and the fact that  $\oF$ is $h^{1+\delta}$-insensitive.
\hfill\qed
\rem{}{ It is easy show that $\frac{h(x)}{x}\searrow0$ for $h(x)$ in the proof of Lemma \hyperref[]{\ref{h-insensitive for F in L}}. 
}

\noindent
\textbf{Proof of Theorem \ref{uniform h-insensitive for convolution in L}.}
 Assume that $\oF_i$ is $h_{i}$-insensitive, where $h_{i}(x)=h(x;F_i)$ is given in Lemma \ref{h-insensitive for F in L}, $1\le i\le n$.
Let $h(x):=h(x; F_1, \cdots, F_n)=\min\{ h_{i}(x), 1\le i\le n\}=o(x)$. Then all $\oF_i$'s are $h$-insensitive and $h(x)\le c h(\frac{x}{c})$,  $c\ge 1,$ by
Lemma \ref{h-insensitive for F in L}. The uniform asymptotic relation \hyperref[]{\eqref{uniform h-insensitive}}, which is essentially the case of $n=2$ in  proof, will be proved by induction. It is obviously true for $n=1$ by Lemma \ref{h-insensitive for F in L}. Since distribution functions are nondecreasing, \hyperref[]{\eqref{uniform h-insensitive}} is equivalent to
\begin{eqnarray}\label{uniform h-insensitive n=2 x+h(x)}
 \lim_{x\toi} \inf_{a(x)\le c_i\le b(x), 1\le i\le n}
         \frac{P\big(\sum_{i=1}^nc_iX_i>x+ h(x)\big)}{P\big(\sum_{i=1}^nc_iX_i>x \big)}\ge1,
\end{eqnarray}
and
\begin{eqnarray}\label{uniform h-insensitive n=2 x-h(x)}
 \lim_{x\toi} \sup_{a(x)\le c_i\le b(x), 1\le i\le n}
          \frac{P\big(\sum_{i=1}^nc_iX_i>x- h(x)\big)}{P\big(\sum_{i=1}^nc_iX_i>x \big)}\le1.
\end{eqnarray}
Write $A+B+C$ for the union of disjoint sets $A,B,C$. The fact that $\big\{\sum_{i=1}^{n}c_iX_i>x\pm h(x)\big\}
=\big\{\sum_{i=1}^{n}  c_{i}X_{i}>x+h(x), c_{n}X_{n}\le \frac{x+h(x)}{2}\big\}
 +\big\{\sum_{i=1}^{n} c_{i}X_{i}>x+h(x),\sum_{i=1}^{n-1}c_{i}X_{i}\le \frac{x+h(x)}{2}\big\}
 +\big\{\sum_{i=1}^{n-1}c_{i}X_{i}> \frac{x+h(x)}{2}, c_{n}X_{n}> \frac{x+h(x)}{2}\big\}$
 and   independence of $X_i$'s yield
\begin{eqnarray}\label{lower bound in Theorem 1}
P\Big(\sum_{i=1}^{n}c_iX_i>x+ h(x)\Big)&\ge&\int_{-\infty}^{x/2} P\Big(\sum_{i=1}^{n-1}c_iX_i>x+h(x)-t\Big) d P(c_{n}X_{n}\le t) \nonumber\\
&&\quad
    +\int_{-\infty}^{x/2}P\big(c_{n}X_{n}>x+h(x)-t\big)d P\Big(\sum_{i=1}^{n-1}c_iX_i\le t\Big) \nonumber\\
&&\quad +P\Big(\sum_{i=1}^{n-1}c_iX_i>\frac{x+h(x)}{2}\Big)P\Big(c_{n}X_{n}>\frac{x+h(x)}{2}\Big).
\end{eqnarray}
The induction assumption with  $b(x)$ replaced by $2b(x)$ implies that 
\begin{eqnarray}\label{decomposition.part.3}
&&P\Big(\sum_{i=1}^{n-1}c_iX_i>\frac{x+h(x)}{2}\Big)P\Big(c_{n}X_{n}>\frac{x+h(x)}{2}\Big)\nonumber\\
&=&P\Big(\sum_{i=1}^{n-1}2c_iX_i>x+h(x)\Big)P\Big(2c_{n}X_{n}>x+h(x)\Big)\nonumber\\
&\sim& P\Big(\sum_{i=1}^{n-1}2c_iX_i>x\Big)P\Big(2c_{n}X_{n}>x\Big)=P\Big(\sum_{i=1}^{n-1}c_iX_i>\frac{x}{2}\Big)P\Big(c_{n}X_{n}>\frac{x}{2}\Big)
\end{eqnarray}
holds uniformly for $a(x)\le  c_i\le  b(x), 1\le i\le n$.\\
Use monotonicity of any distribution function and the inequality $h(x)\le 2h(\frac{x}{2})$ to obtain
\begin{eqnarray}\label{decomposition.part.1}
1\ge\inf_{t\le x/2}\frac{\oF(x+h(x)-t)}{\oF(x-t)}
\ge \inf_{t\le x/2}\frac{\oF\big(x-t+2h(\frac{x}{2})\big)}{\oF(x-t)}
\ge\inf_{u=x-t\ge x/2}\frac{\oF(u+2h(u))}{\oF(u)}\sim 1
\end{eqnarray}
provided $\overline{F}$ is $h$-insensitive. It follows from the induction assumption and Lemma \ref{h-insensitive for F in L} that the tail distribution of $\sum_{i=1}^{n-1}c_iX_i$ and the tail distribution of $c_{n}X_{n}$ are $h$-insensitive.
The asymptotic relation \hyperref[]{\eqref{decomposition.part.3}} and the inequality \hyperref[]{\eqref{lower bound in Theorem 1}} imply
\begin{eqnarray*}
&&P\big(\sum_{i=1}^{n}c_iX_i>x+ h(x)\big)\\
&\ge&\bigg(\int_{-\infty}^{x/2} P\Big(\sum_{i=1}^{n-1}c_iX_i>x-t\Big) d P(c_{n}X_{n}\le t)
    +\int_{-\infty}^{x/2}P\big(c_{n}X_{n}>x-t\big)d P\Big(\sum_{i=1}^{n-1}c_iX_i\le t\Big) \nonumber\\
&&\quad +P\Big(\sum_{i=1}^{n-1}c_iX_i>\frac{x}{2}\Big)P\Big(c_{n}X_{n}>\frac{x}{2}\Big)\bigg)(1+o(1))\\
&=& (1+o(1))P\Big(\sum_{i=1}^{n}c_iX_i>x\Big),
\end{eqnarray*}
where the term $o(1)$ goes to 0  uniformly for $a(x) \le c_{i}\le b(x)$, $1\le i\le n$. This complete the proof of  \hyperref[]{\eqref{uniform h-insensitive
n=2 x+h(x)}}.\\
The other uniform asymptotic relation \hyperref[]{\eqref{uniform h-insensitive n=2 x-h(x)}} can be obtained by substituting  $+h(x)$,
$+2h(\frac{x}{2})$, $\ge$, $\inf$ with $-h(x)$, $-2h(\frac{x}{2})$, $\le$, $\sup$, respectively, in the proof of \hyperref[]{\eqref{uniform h-insensitive
n=2 x+h(x)}}.
\hfill\qed


\bigskip
\noindent
{\bf Proof of Theorem \hyperref[]{\ref{2nd theorem}}.} The idea is from the proof of Theorem 2.1 of Chen et al. \hyperref[]{\cite{CNY11}}.
Let  $\big\{\Omega_K=\{X_i\ge 0\mathrm{\ for\ all\ }i\in K, X_j< 0\mathrm{\ for\ all\ }j\in \{1,\cdots,n\}\backslash K\}, K\subseteq\{1,\cdots,n\}\big\}$ be a finite partition of the whole space $\Omega$. Obviously, $P\big(\sum_{i=1}^n c_iX_i>x,\Omega_K\big)$ is not less than
\begin{eqnarray}\label{2nd theorem 1st equation}
 & & P\Big(\sum_{i\in K}c_iX_i>x+h(x),\sum_{j\notin K}c_jX_j>-h(x) ,\Omega_K\Big)\nonumber\\
&=&P\Big(\sum_{i=1}^nc_iX_i^+>x+h(x),\Omega_K\Big)-P\Big(\sum_{i\in K}c_iX_i>x+h(x),\sum_{j\notin K}c_jX_j\le -h(x) ,\Omega_K\Big),
\end{eqnarray}
where, due to the independence of $X_i$'s, the second term   equals
\be
 P\Big(\sum_{i\in K}c_iX_i>x+h(x),\bigcap_{i\in K}\{X_i\ge0\}\Big)P\Big(\sum_{j\notin K}c_j (-X_j) \ge h(x), \bigcap_{j\notin K}\{X_j<0\}\Big).
\ee
and it is at most $P\big(\sum_{i=1}^n c_iX_i^+>x+h(x)\big)P\big(\sum_{j=1}^n c_jX_j^-\ge h(x)\big)$, where $x^-=\max\{-x,0\}$. Note that
$\{\sum_{j=1}^n c_jX_j^-\ge h(x)\}\subseteq\bigcup_{j=1}^n\{c_jX_j^-\ge\frac{h(x)}{n}\}=\bigcup_{j=1}^n\{c_jX_j\le-\frac{h(x)}{n}\}$,
whose probability is at most $\sum_{j=1}^n  P\Big(X_j\le -\frac{h(x)}{nb(x)}\Big)=o(1)$ provided $b(x)=o(h(x))$. Therefore, uniformly for $0<a\le c_i\le b(x)$, $1\le i\le n$, the second term in \hyperref[]{\eqref{2nd theorem 1st equation}} is $o\big(P\big(\sum_{i=1}^n c_iX_i^+>x+h(x)\big)\big)$
and
\begin{eqnarray*}
P\Big(\sum_{i=1}^n c_iX_i>x,\Omega_K\Big) \ge P\Big(\sum_{i=1}^n c_iX_i^+>x+h(x),\Omega_K\Big) +o\Big(P\Big(\sum_{i=1}^n c_iX_i^+>x+h(x)\Big)\Big).
\end{eqnarray*}
Sum it over all $K$'s to get
\be
P\Big(\sum_{i=1}^n c_iX_i>x\Big) \ge P\Big(\sum_{i=1}^n c_iX_i^+>x+h(x)\Big) +o\Big(P\Big(\sum_{i=1}^n c_iX_i^+>x+h(x)\Big)\Big).
\ee
Clearly, $X_i^+\sim F_i^+(x)=F_i(x)I_{\{x\ge0\}}\in\mathcal{L}, 1\le i\le n$.  Choose $h(x)$ such that \hyperref[]{\eqref{uniform h-insensitive}} holds with $F_i$ substituted by $F_i^+$.
 The desired result follows from Theorem \ref{uniform h-insensitive for convolution in L} and the simple fact that $\sum_{i=1}^{n} c_iX_i\le \max_{1\le k\le n}\sum_{i=1}^k c_iX_i \le \sum_{i=1}^{n} c_iX_i^+$.
 \hfill\qed

\noindent
{\bf Proof of Corollary \hyperref[]{\ref{S in R}}.}  Recall that $\overline{F}\in\mathcal{S}$ if
$\overline{F^+}\in\mathcal{S}$, i.e., $\overline{F^+*F^+}(x)\sim 2\overline{F^+}(x)$ for $F^+(x)=F(x)I_{\{x\ge0\}}$.
Clearly, $F\in\mathcal{L}$ iff $F^+\in\mathcal{L}$. If $F^+\in\mathcal{S}$, the fact that $\mathcal{S}\subset\mathcal{L}$   implies $F\in\mathcal{L}$.
Then it is equivalent to show that $\overline{F^+*F^+}(x)\sim 2\overline{F^+}(x)$ iff $\overline{F*F}(x)\sim 2\overline{F}(x)$, i.e. $\overline{F^+*F^+}(x)\sim \overline{F*F}(x)$ since $\overline{F^+}(x)=\overline{F}(x)$ for all $x>0$. It is obviously true by Theorem $\hyperref[]{\ref{2nd theorem}}$.
\hfill\qed

\bigskip
The next two lemma can be easily checked from the definition of the class $\mathcal{C}$.
\lem{}{\label{cX is h-insensitive for any 0<c<pb(x)=o(x)}
If $X$ follows distribution $F\in\mathcal{C}$, then $\oF(x)$ is $h$-insensitive provided $h(x)=o(x)$ and it holds that, uniformly for $0<c<b(x)=o(x)$,
\be
P(c X>x\pm h(x))\sim P(c X>x).
\ee
}
\lem{}{\label{independent for TAI implication}
If $X_i\sim F_i\in\mathcal{C}, 1\le i\le n$, are pQSAI random variables,  it holds that, uniformly for $0<c<b(x)=o(x)$,
\be
 P\Big( c_jX_j>\frac{x}{n},\max_{1\le k\ne j\le n}|c_kX_k|> b(x)\ln\Big(\frac{x}{b(x)}\Big) \Big)=o\big(P(c_jX_j>x)\big)
\ee
and consequently
\be
 P\Big(\bigcup_{j=1}^n\Big\{c_jX_j>\frac{x}{n},\max_{1\le k\ne j\le n}|c_kX_k|> b(x)\ln\Big(\frac{x}{b(x)}\Big)\Big\}
 \Big)=o\Big(\sum_{j=1}^nP(c_jX_j>x)\Big).
\ee
}
\noindent
{\bf Proof of Theorem \ref{3rd theorem}.}   Let $h(x)=b(x)\ln\big(\frac{x}{b(x)}\big)$. The proof is similar to that of Theorem \ref{4th theorem} and is
omitted. \hfill\qed

\noindent
{\bf Proof of Corollary \ref{Cor to 3rd Theorem}.}  Partition the range of the weights as
$\{(c_1,\cdots,c_n): 0\le c_i\le b(x), 1\le i\le n, \min_{i=1}^{n}c_i>0\} =\bigcup_{K\subset\{1,\dots,n\}}\{(c_1,\cdots,c_n): 0\le c_i\le b(x), i\in K, 0< c_i\le b(x), i\notin K \}$. The desired result follows from Theorem \ref{3rd theorem}.
\hfill\qed

\lem{}{\label{TAI implication}
If $X_i\sim F_i\in\mathcal{D}, 1\le i\le n$,  are pSQAI random variables, $h(x)=o(x)$ and $h(x)\nearrow\infty$, it holds that, uniformly for $0<a<c_i<b(x)=o(h(x)), 1\le i\le n$,
\be
 P\Big( c_jX_j>\frac{x}{n},\max_{1\le k\ne j\le n}|c_kX_k|> h(x) \Big)=o\big(P(c_jX_j>x)\big)
\ee
and consequently
\be
 P\Big(\bigcup_{j=1}^n\big\{c_jX_j>\frac{x}{n},\max_{1\le k\ne j\le n}|c_kX_k|> h(x)\big\} \Big)=o\Big(\sum_{j=1}^nP(c_jX_j>x)\Big).
\ee
}
\noindent
{\bf Proof.} The results follow from the fact that $F_i\in\mathcal{D}$ and $b(x)=o(h(x))$, the pSQAI property of $X_i$'s  and the elementary probability
inequality $P(A\cap\cup_{i=1}^nB_i)\le \sum_{i=1}^n P(AB_i)$.
\hfill \qed

If  $X_i$  is large,  the pSQAI property of $X_j$'s implies that other $X_j$'s are relatively close to 0 and negligible compared with $X_i$.
If $\sum_{i=1}^n c_iX_i>x$,  there should be exactly one $c_iX_i$ greater than $\frac{x}{n}$ and consequently Lemma \hyperref[]{\ref{TAI implication}} implies
\be
P\Big(\sum_{i=1}^n c_iX_i>x \Big)
\sim \sum_{j=1}^n P\Big(\sum_{i=1}^n c_iX_i>x, c_jX_j>\frac{x}{n},\max_{1\le k\ne j\le n}|c_kX_k|\le h(x)  \Big).
\ee
It gives the idea of the proof of  Theorem \ref{4th theorem},
which is simpler and more straightforward than the proof of Lemma 2.1 of Liu et al. \hyperref[]{\cite{LGW12}} and Theorem 2.1 of Li \hyperref[]{\cite{L13}}.

\noindent
{\bf Proof of Theorem \ref{4th theorem}.} All asymptotic relations hold uniformly for $a(x)\le c_i\le b(x), 1\le i\le n,$ in the proof. By
Lemma \hyperref[]{\ref{h-insensitive for F in L}},  there exists a positive nondecreasing function $h(x):=h(x,a; F_1,\cdots,F_n)$ satisfying $h(x)\nearrow\infty$ and
$h(x)=o(x)$ such that
\hyperref[]{\eqref{Theorem3.1.n=1}} holds for $F=F_i,1\le i\le n$, respectively. Choose $b(x)=o(h(x))$ and $b(x)\nearrow\infty$.
Note that
\be
\Big\{\sum_{i=1}^n c_iX_i>x \Big\}&=&\bigcup_{j=1}^n \Big\{\sum_{i=1}^n c_iX_i>x ,c_jX_j>\frac{x}{n}\Big\}\\
 &=&\bigcup_{j=1}^nA_j \bigcup \Big\{\sum_{i=1}^n c_iX_i>x, \bigcup_{j=1}^n\big\{c_jX_j>\frac{x}{n},\max_{1\le k\ne j\le n}|c_kX_k|> h(x)\big\} \Big\},
\ee
where $A_j=\big\{\sum_{i=1}^n c_iX_i>x, c_jX_j>\frac{x}{n},\max_{1\le k\ne j\le n}|c_kX_k|\le h(x)\big\}, 1\le j\le n$, are mutually exclusive events
provided $\frac{x}{n}>h(x)$.
The elementary probability inequality $P(A)\le P(A\cup B)\le P(A)+P(B)$ and Lemma \hyperref[]{\ref{TAI implication}} lead to
\begin{eqnarray}\label{Decompostion of probability of large sum}
P\Big(\sum_{i=1}^n c_iX_i>x \Big)=\sum_{j=1}^n P(A_j)  +o\Big(\sum_{j=1}^nP(c_jX_j>x)\Big).
\end{eqnarray}
Lemma \hyperref[]{\ref{h-insensitive for F in L}} and the fact that $c_jX_j$ is at least $x-(n-1)h(x)$ on $A_j$
lead to
\begin{eqnarray*}
 P(A_j)\le  P\big(c_jX_j>x-(n-1)h(x)\big) = P(c_jX_j>x)+o\big(P(c_jX_j>x)\big), \quad 1\le j\le n.
\end{eqnarray*}
Since $\max_{1\le k\ne j\le n}|c_kX_k|\le h(x)$ on $A_j$, $c_j X_j > x + (n -1)h(x)$ implies $\sum_{i=1}^n c_iX_i>x$ on $A_j$ for any $1\le j\le n$.
It follows from Lemma \hyperref[]{\ref{h-insensitive for F in L}}  and \hyperref[]{\ref{TAI implication}}  that
\begin{eqnarray*}
 P(A_j) &\ge &  P\big(c_jX_j>x+(n-1)h(x),\max_{1\le k\ne j\le n}|c_kX_k|\le h(x)\big)\\
 &=&P(c_jX_j>x+(n-1)h(x))- P\big(c_jX_j>x+(n-1)h(x),\max_{1\le k\ne j\le n}|c_kX_k|> h(x)\big)\\
&=& P(c_jX_j>x)+o\big(P(c_jX_j>x)\big),\quad 1\le j\le n.
\end{eqnarray*}
Therefore, \hyperref[]{\eqref{Decompostion of probability of large sum}} can be written as
\begin{eqnarray} \label{1st limit in Theorem 4}
 P\Big(\sum_{i=1}^{n} c_iX_i>x \Big)\sim \sum_{i=1}^{n} P(c_iX_i>x).
\end{eqnarray}
In the exactly same way, it can be proved that
\begin{eqnarray}\label{2nd limit in Theorem 4}
 P\Big(\sum_{i=1}^{n} c_iX_i^+>x \Big)\sim \sum_{i=1}^{n} P(c_iX_i^+>x)=\sum_{i=1}^{n} P(c_iX_i>x).
\end{eqnarray}
Note that $\sum_{i=1}^{n} c_iX_i\le \max_{1\le k\le n}\sum_{i=1}^k c_iX_i \le \sum_{i=1}^{n} c_iX_i^+$. The desired results follow from the uniform
asymptotic relation \hyperref[]{\eqref{1st limit in Theorem 4}} and  \hyperref[]{\eqref{2nd limit in Theorem 4}}.
 \hfill\qed

\rem{}{The proof of Theorem \ref{4th theorem} also leads to Corollary \ref{liminf for L}.
}

\bigskip
\noindent
\textbf{Acknowledgments}\\
The author would like to thank the anonymous referees for their comments  and help in improving the paper.


\end{document}